\documentclass
{book}
\usepackage[editedvolume,numreferences]{crckapb} %
\usepackage{amsmath,amssymb,amscd,amsthm,latexsym}
\newtheorem{theorem}{Theorem}[section]

\newtheorem*{rmk}{Remark}
\newtheorem*{srmk}{Remarks}

\begin{opening}
\title{A TREE-ARROWING GRAPH
}
\author{S. Shelah${}^1$}
\institute{ Hebrew
University,\\ Jerusalem, Israel.}
\author{E. C. Milner${}^2$}
\institute{University of Calgary,\\ Calgary, Canada.}
\runningauthor{S. Shelah and E. C. Milner}
\end{opening}
\begin{document}
\footnotetext[1]{Paper Sh 578 in Shelah's publication list.  
Research supported by ``The Israel Science Foundation'' administered 
by The Israel Academy of Sciences and Humanities.}

\footnotetext[2]{Research supported by NSERC grant \#69-0982.} 

\addtocounter{footnote}{2}

{\it Dedicated to the memory of Eric Milner}

\bigskip
\begin{abstract}
We answer a variant of a question of R\" odl and Voigt by showing
that, for a given infinite cardinal $\lambda$, there is a graph $G$ of
cardinality $\kappa=(2^\lambda)^+$ such that for any colouring of the
edges of $G$ with $\lambda$ colours, there is an induced copy of the
$\kappa$-tree in $G$ in the set theoretic sense with all edges having
the same colour.

\end{abstract}

\section{Introduction}

 ${\cal   G}=(V,E)$ is a graph with {\it vertex set $V$} 
and {\it edge set $E$},  where $E \subseteq  [V]^2$.  
The graph ${\cal   H}=(W,F)$ is a {\it subgraph\/} of 
${\cal   G}$ if $W \subseteq  V$ and $F \subseteq   E$, 
it is an {\it induced subgraph \/} if $F=E\cap [W]^2$.  
If $\lambda$ is a cardinal, the partition relation
\begin{equation} {\cal   G} \rightarrow ({\cal   H})^2_ \lambda 
\label{part}, \end{equation}
means that if $c:E \rightarrow \lambda$ is any colouring of the 
edges of ${\cal   G}$ with $\lambda$ colours, then there is an induced copy of 
${\cal   H}$ in ${\cal   G}$ in which all the edges have the same colour.
There is a related notion ${\cal   G} \rightarrow ({\cal   H})^1_
\lambda$, 
for vertex colourings of graphs.  However, there is an essential difference
since, for any given graph ${\cal   H}$ and any $\lambda$, there 
is some ${\cal   G}$ such that  
${\cal   G} \rightarrow ({\cal   H})^1_ \lambda$
holds.  This is not true for edge-colourings; 
Hajnal and Komjath \cite{HK} proved the consistency of a negative
answer, 
and Shelah \cite{Sh:289} proved that a positive answer is 
also consistent.  It is therefore of some interest to have 
instances of graphs ${\cal   H}$ such that (\ref{part}) holds 
for some ${\cal   G}$, and then, of course, one can ask for the 
smallest such ${\cal   G}$.

\medskip

R\"{o}dl and Voigt \cite{RV} (see also \cite{KM}) proved a 
result of this kind by showing that
for any infinite cardinal $\lambda$ and a suitably large 
$\kappa$, there is a graph  ${\cal   G}_ \kappa$ of cardinality 
$\kappa$ such that 
\begin{equation}{\cal   G}_ \kappa \rightarrow 
({\cal   T_ \kappa})^2_ \lambda \label{tree}\end{equation}
holds, where ${\cal   T}_ \kappa$ is the tree in which every 
vertex has degree $\kappa$ (see below).  More precisely, 
`suitably large'  means that the ordinary partition relation
\[ {\rm cf}(\kappa) \rightarrow (\omega)^3_ \lambda\]
holds so that, by \cite{ER}, $\kappa \geq (2^{2^\lambda})^+$; 
in fact, they showed in this case that the ubiquitous 
{\it shift-graph\/} on $\kappa$ works.
R\"{o}dl and Voigt \cite{RV} then asked, what is the smallest 
cardinal $\kappa$
such that (\ref{tree}) holds?  It is easily seen that (\ref{tree}) 
is false 
if $\kappa\leq 2^\lambda$, and they conjectured that it holds 
(for some suitable graph ${\cal   G}_ \kappa$) if $\kappa=(2^\lambda)^+$.
In this paper we prove that
(\ref{tree})  holds with ${\cal   T}_ \kappa$ replaced by 
${\cal   T}(\kappa)$, a related graph which we call 
{\it the transitive $\kappa$-tree\/} defined in the next section.

\section{Preliminaries}
For an infinite cardinal $\kappa$ we denote by ${}^{<\omega}\kappa$
 the set of all increasing finite sequences of ordinals in $\kappa$.  
The {\it length\/} 
of an element $s=\langle s_0,\ldots,s_{n-1}\rangle\in 
{}^{<\omega}\kappa$
 is denoted by 
$\ell n(s)=n$.  Also, we define 
$$\max(s)=\left \{\begin{array}{ll} -1 &\mbox{ if $s=\langle 
\rangle$, the empty sequence,}\\
s_{\ell n(s)-1} & \mbox{ if $\ell n(s)>0$.} 
\end{array}\right.$$
If $s=\langle s_0,\ldots,s_{n-1}\rangle$ and 
$t=\langle t_0,\ldots,t_{m-1}\rangle$ are two elements of 
${}^{<\omega}\kappa$, we write $s\lhd t$ to denote the fact that $s$ is a
proper initial segment of $t$, that is $n<m$ and $s_i=t_i$ for $i<n$, and 
in this case we write $s=t|n$.  We also write $s=t_*$ if 
$m=n+1$ and $s\lhd t$.  If $s,t$ are distinct and 
$\lhd$-incomparable we write $s\perp t$.  
The {\it $\kappa$-tree of height $\omega$} is the graph ${\cal   T}_\kappa$ 
on ${}^{<\omega}\kappa$
with edge set 
$$E_ \kappa =\{\{s,t\}: s,t\in {}^{<\omega}\kappa \wedge s=t_*\}.$$
We shall also consider a related graph, {\it the transitive 
$\kappa$-tree of height $\omega$}, which is  the graph ${\cal   T}(\kappa)$ 
on ${}^{<\omega}\kappa$
with edge set 
\[F_ \kappa=\{\{s,t\}:s,t\in {}^{<\omega}\kappa \wedge s\lhd t\}.\]

\medskip

We shall prove the following theorem. 
\begin{theorem} Let $\lambda$ be an infinite cardinal, 
and let $\kappa=(2^\lambda)^+$.  Then there is a 
graph $G_ \kappa$ of cardinality $\kappa$ such that
\[G_ \kappa\rightarrow ({\cal  T})^2_ \lambda,
\]
where ${\cal   T}$ is 
${\cal   T}(\kappa)$.\label{main}\end{theorem}

\begin{rmk}
Instead of $\kappa=(2^\lambda)^+$, 
it is enough that $\kappa$ be any regular cardinal 
such that $|\alpha|^\lambda<\kappa$ holds for all 
$\alpha<\kappa$.  The same proof works.
\end{rmk}

The construction of a suitable ${\cal   G}_ \kappa$ 
depends upon the following (slightly weaker version of a) 
theorem of Shelah \cite{Sh:413} (or more \cite[3.5]{Sh:572}):
\begin{quote} $(\bullet)$ {\it  Let
$\lambda$ be an infinite cardinal, $\kappa=(2^\lambda)^+$,
$S=\{\alpha<\kappa: {\rm cf}(\alpha)=\lambda^+\}$.  Then
there are a sequence $\overline{C}=\langle  C_ \delta:\delta\in S \rangle$ 
and a sequence $\overline{h^*}=\langle  
h^*_ \delta:\delta\in S \rangle$ such that
$C_ \delta$   is a club in $\delta$ having 
order type $\lambda^+$, $h^*_ \delta:C_ \delta \rightarrow 2$
and such that, for any club $K$ in $\kappa$, there is a 
stationary subset $B_K$ of $S\cap K$
such that for each $\delta\in B_K$ and each $i<2$,
$\min(C_ \delta)\in K$ and the set
$$D_K{(\delta,i)}=\{\alpha\in C_ \delta\cap K:h^*_ \delta(\alpha)
=i \wedge \min(C_ \delta\,\backslash\,(\alpha+1))\in K\}$$
is cofinal in $\delta$.}
\end{quote} 

\begin{srmk}
1. The result is also true if $2$, 
the range of each $h^*_ \delta$, is replaced by $\lambda$; 
also, if $\kappa=\lambda^{++}$,
we can also require that
$D_{K}(\delta,i)$ be a stationary subset of $\delta$ for each $\delta\in B_K$
and $i<\lambda$ (see \cite{Sh:572}).

2. If $2^\lambda>\lambda^+$, then the following stronger 
assertion is true
(see Shelah \cite{Sh:365}): $(\bullet\bullet)${\it There 
is a sequence $\bar{C}=\langle  C_ \delta:\delta\in S \rangle$ such that
$C_ \delta$   is a club in $\delta$ having order type 
$\lambda^+$ and, for any club $K$ in $\kappa$
and any stationary  subset $S'\subseteq  S$, 
there is a stationary 
subset $B_K \subseteq  S'\cap K$
such that $C_ \delta \subseteq   K$ for each $\delta\in B_K$.}  
Using this result instead of $(\bullet )$, the proof of 
Theorem~\ref{main} for the case when $2^\lambda>\lambda^+$ 
may be slightly simplified.
\end{srmk}

\medskip

 We will prove that Theorem~\ref{main} holds with the graph 
$G_ \kappa=(\kappa,{\cal   E})$, where
\[ {\cal   E}=\{\{\alpha,\beta\}: \beta\in S \wedge 
\min(C_ \beta) <\alpha<\beta
\wedge h^*_ \beta(\sup(\alpha\cap C_ \beta))=0\},\]
and the $C_ \beta$ and $h^*_ \beta$ are as described in ($\bullet$).

\section{The case ${\cal   T}={\cal   T}(\kappa)$}
We prove the result for the case of the transitive tree $T(\kappa)$.
\medskip

\noindent{\it Proof: } 
Let $c:{\cal   E} \rightarrow \lambda$ be any $\lambda$-colouring 
of the edges of $G_ \kappa$.  For each $\zeta\in \lambda$ consider 
the following two-person game ${\cal  G}_ \zeta$.  The game has 
$\omega$ moves.  At the $n$-th stage the first player $P_1$ 
chooses ordinals $\alpha_n, \beta_n$, and then the second 
player $P_2$ chooses two ordinals $\gamma_n, \delta_n$ so that 
\begin{equation} \alpha_n< \beta_n < \gamma_n < \delta_n<\kappa, 
\label{G1}\end{equation}
\begin{equation} \delta_m<\alpha_n \quad (m<n). \end{equation}
The player $P_2$ is declared the winner in a play of the game 
if he succeeds in choosing the $\gamma_n$ so that
\begin{equation} \{\gamma_m,\gamma_n\}\in {\cal   E}, \qquad
c(\{\gamma_m,\gamma_n\})= \zeta \quad (m<n<\omega), 
\label{W2} \end{equation}
and
\begin{equation} \{\xi,\gamma_n\}\notin {\cal   E}  
\mbox{\ \ for \  $\xi\in (\alpha_m,\beta_m)$ and  
$m\leq n<\omega$.}\label{W3} \end{equation}
(As usual, $(\alpha,\beta)$  denotes the open interval
$\{\xi:\alpha<\xi<\beta\}$ and $[\alpha,\beta]$ is the 
corresponding closed interval.)

\medskip
The proof of the theorem depends upon
 the following two facts:

\noindent {\bf Fact A:}\  For some $\zeta< \lambda$, 
$P_2$ has a winning strategy for the game ${\cal G}_ \zeta$.
\medskip

\noindent{\bf Fact B:}\ If $P_2$ can win ${\cal  G}_\zeta$, then the graph $G_ \kappa$ contains an induced copy of ${\cal  T}(\kappa)$ with all edges coloured $\zeta$.

\begin{proof}
[Proof of Fact B]
We assume that $\zeta< \lambda$ and that the second player $P_2$ has a winning strategy $\sigma_ \zeta$ for the game ${\cal G}_\zeta$.  We shall define ordinals $\alpha_s, \beta_s, \gamma_s, \delta_s$ for $s$ a vertex of ${\cal    T}(\kappa)$ so that the 
following conditions are satisfied:

\noindent\newline (a) For each $s$  the sequence
\[\langle  (\alpha_{s|i}, \beta_{s|i},\gamma_{s|i},\delta_{s|i}): i<\ell n(s) \rangle\]
consists of the first $2\ell n(s)$ moves in a proper play of the game ${\cal G}_\zeta$ in which $P_2$ uses the winning strategy $\sigma_ \zeta$.

\noindent\newline (b)  $\gamma_s\neq \gamma_t$ if $s\neq t$.

\noindent\newline (c) If $s\perp t$, then $\{\gamma_s, \gamma_t\}\notin {\cal   E}$.

\medskip

 Since (\ref{W2}) holds, these conditions imply that the map $s \mapsto \gamma_s$ is an embedding of the tree ${\cal   T}(\kappa)$ into the graph $G_ \kappa$ and all the edges of the image have colour $\zeta$.

\medskip
In fact, we shall choose the $\alpha_s, \beta_s, \gamma_s, \delta_s$ so that (a) holds and so that the following condition is satisfied:

\noindent\newline (d) For any vertices $s, t$ of ${\cal   T}(\kappa)$, if
 $s\perp t$, then 
\begin{tabbing}
EITHER \={\it (i)} \=$\quad [\gamma_s, \delta_s] \subset \bigcup_{i\leq\ell n(t)}(\alpha_{t|i}, \beta_{t|i})$,\\
OR \>{\it (ii)} \>$\quad [\gamma_t, \delta_t] \subset \bigcup_{i\leq\ell n(s)}(\alpha_{s|i}, \beta_{s|i})$.
\end{tabbing}

\noindent  The conditions (a) and (d), and the fact that $P_2$ is using the winning strategy $\sigma_\zeta$,  ensure that (b) and (c) also hold.

\medskip
We define $\alpha_s, \beta_s, \gamma_s, \delta_s$ by induction on $\max(s)$. 
Let $\alpha_{\langle \rangle}=0$, $\beta_{\langle \rangle}=1$, and then let 
$(\gamma_{\langle \rangle}, \delta_{\langle \rangle})$ be $P_2$'s response in the game ${\cal  G}_ \zeta$ using his winning strategy $\sigma_ \zeta$.
Now let $0\leq\xi< \kappa$, and suppose that we have suitably defined
 $\alpha_s, \beta_s, \gamma_s, \delta_s$ for all vertices $s$ of
${\cal   T}(\kappa)$ such that
$\max(s)<\xi$.  We need to define these when $\max(s)=\xi$.

\medskip  Let $\langle t_i:i<\theta(\xi) \rangle$ be an enumeration of all 
the nodes $s$ of ${\cal T}(\kappa)$ with $\max(s)=\xi$.  Then $1\leq \theta(\xi)\leq
2^\lambda< \kappa$.  Now inductively choose the $\alpha_{t_i}$, $\beta_{t_i}$, $\gamma_{t_i}$, $\delta_{t_i}$ for $i<\theta(\xi)$ so that
\[ \alpha_{t_i}= \delta_{(t_i)_*}+1,\]
and if $i=0$, $\beta_{t_i} = \alpha_{i_0}+1$ and if $i>0$
\[\beta_{t_i}= \sup\{ \delta_s +2: \max(s)<\xi \mbox{ or  }
s=t_j \mbox{ for some } j<i\}.\]
The corresponding pairs $(\gamma_{t_i}, \delta_{t_i})$ are determined by the strategy $\sigma_ \zeta$.  With these choices it is easily seen that (a) continues to hold; we have to check that (d) also holds when $s\perp t$ and $\max(s)=\xi$ or $\max(t)=\xi
$.

\medskip
If $\max(s)=\max(t)=\xi$, then $s=t_i$ and $t=t_j$, where say $i<j$.  Then 
\[ \alpha_t=\delta_{t_*}+1< \beta_s<\gamma_s<\delta_s<\beta_t,\]
and so (d)(i) holds.  

\medskip 
Suppose $\max(s)<\xi=\max(t)$.  Then by the induction hypothesis, either (i) or (ii) of (d) holds when we replace $t$ by $t_*$.
Suppose first that (d)(i) holds. Then for some $m\leq\ell n(t_*)$ we have that
\[\alpha_{t_*|m}<\gamma_s<\delta_s<\beta_{t_*|m}.\]
  It follows that (d)(i) also holds for $s$ and $t$ since $t|m=t_*|m$.  Now suppose that 
(d)(ii) holds so that, for some $m\leq\ell n(s)$,
\[ \alpha_{s|m}<\gamma_{t_*}<\delta_{t_*}<\beta_{s|m}.\]
Then, by the definitions of $\alpha_t$ and $\beta_t$, it follows that
\[\alpha_t=\delta_{t_*}+1\leq\beta_s<\gamma_s<\delta_s<\beta_t,\]
so that again (d)(i) holds for $s$ and $t$.  Similarly, if $\max(t)<\xi=
\max(s)$.\ \  
\end{proof}

\begin{proof}[Proof of Fact A]  
We have to show that $P_2$ wins the game ${\cal   G}_\zeta$ for some 
$\zeta<\lambda$.  Suppose for a contradiction that this is false.  
Since the games are open and hence determined, it follows that $P_1$ 
has a winning strategy, say $\tau_\zeta$, for the game
${\cal   G}_\zeta$ for every $\zeta<\lambda$. 
\medskip

For convenience we write $c(\{\alpha,\beta\})=-1$ if 
 $\{\alpha,\beta\}\notin {\cal   E}$, so that $c$ is defined 
on all pairs  $\{\alpha,\beta\}\in[\kappa]^2$.
For each bounded subset $X \subseteq  \kappa$ define an equivalence 
relation $e_X$ on $S \,\backslash\, (\sup (X )+1)$ so that 
$\beta\, e_X \,\gamma$ holds if and only if 
\begin{enumerate}
\item[]\begin{enumerate}\renewcommand{\theenumii}{%
\roman{enumii}}
\item 
$\beta,\gamma\in S$ 
and $\sup(X)<\beta,\gamma<\kappa$; 
\item 
$c(\{\alpha,\beta\})=
c(\{\alpha,\gamma\})$ for all $\alpha\in X$; 
\item 
$X\cap C_ \beta
=X\cap C_ \gamma$, (iv) for $\alpha\in X$, $\alpha\leq\min(C_ \beta) 
\Leftrightarrow  \alpha\leq\min (C_ \gamma)$,
 $ \mathop{\rm  tp}(\alpha\cap C_ \beta)= 
\mathop{\rm  tp}(\alpha\cap C_ \gamma)$ and 
$h^*_ \beta(\sup(\alpha\cap C_ \beta))=
h^*_ \gamma(\sup(\alpha\cap C_ \gamma))$  
(for $\alpha>\min(C_ \beta)$).  
\end{enumerate}\end{enumerate}
Note that the equivalence relation 
$e_X$ has at most $(\lambda^+)^{|X|}\leq 2^{\lambda|X|}$ classes.  
Also, if $Y \subseteq   X$, then $\beta\, e_X \,\gamma \Rightarrow 
\beta \,e_Y \,\gamma$.

\medskip

Since $\kappa=(2^\lambda)^+$, there is 
a continuous increasing sequence of ordinals \linebreak
$\langle\rho_\eta:\eta<\kappa\rangle$ in $\kappa$ such that the following two conditions hold:

\begin{enumerate}\item[]
\begin{enumerate}
\item[(o)] If $X \subseteq  \rho_\eta$, $|X|\leq \lambda$ and $\rho_\eta<\beta<\kappa$, then there is  some $\gamma\in (\rho_\eta, \rho_{\eta+1})$ such that $\beta e_X \gamma$

\medskip

\item[(oo)] $\rho_\eta$ is closed under $\tau_\zeta$ for all $\zeta<\lambda$.  In other words, if at the $n$-th stage of a play in the game
${\cal   G}_\zeta$, player $P_2$ chooses $\gamma_n<\delta_n<\rho_\eta$, then
$P_1$'s response using $\tau_\zeta$ is to choose $\alpha_{n+1}, \beta_{n+1}$ so that $\delta_n<\alpha_{n+1}<\beta_{n+1}<\rho_\eta$.
\end{enumerate}
\end{enumerate}

\medskip
Since $K=\{\rho_\eta:\eta<\kappa\}$ is a club in $\kappa$, there is some $\delta
\in S$ such that $\min(C_ \delta)\in K$ and, for $\varepsilon\in \{0,1\}$,
$$ A_ \varepsilon=\{ \alpha\in C_ \delta\cap K: h^*_ \delta(\alpha)=\varepsilon \wedge \min(C_ \delta \,\backslash\,(\alpha+1))\in K\}$$
is an unbounded subset of $\delta$.  Let 
$C_ \delta=\{i_ \sigma:\sigma<\lambda^+\}$, where $i_0<i_1<\cdots$\,.

\medskip

We claim that the following assertion holds for some $\zeta< \lambda$.
\begin{quote}

$(*)_\zeta$:  If $X \subseteq  \delta$, $|X|\leq \lambda$, then there are
$\sigma<\lambda^+$ and $\gamma $ such that (a)~$\, \sup(X)<i_ \sigma<\gamma <i_{\sigma+1}$, (b) $i_ \sigma\in A_0$, 
  (c)$\ \gamma\,  e_X\, \delta$, and (d) $c(\gamma ,\delta)=\zeta$.
\end{quote}

 For suppose the claim is false.  Then, for each $\zeta<\lambda$ 
there is a counter-example $X_\zeta$.
Let $X=\bigcup\{X_\zeta:\zeta<\lambda\}$.
Then $X \subseteq   \delta$ and $|X|\leq \lambda$ and so, for some $\alpha\in A_0$, 
$\sup(X)<\alpha<\delta$.  There are $\eta<\kappa$ and $\sigma<\lambda^+$ such that $\alpha= \rho_\eta=i_ \sigma$, and therefore, by the choice of $\rho_{\eta+1}$, there is $\gamma $ 
such that $\rho_\eta<\gamma<\rho_{\eta+1}$ and $\gamma  e_X \delta$.
Since $\alpha=i_ \sigma\in A_0$, $i_{\sigma+1}=\min(C_ \delta \,\backslash\,
(\alpha+1))\in K$.  So $\rho_{\eta+1}\leq i_{\sigma+1}$.  Therefore, $\sup(C_ \delta\cap \gamma )=i_ \sigma$, and since $\alpha=i_ \sigma\in A_0$, we have that $h^*_ \delta(\sup(C_ \delta\cap \gamma ))=0$.  Therefore, $\{\gamma ,\delta\}$ is an edge of $G
$ and there is some $\zeta\in \lambda$ such that $c(\gamma ,\delta)=\zeta$.  But this contradicts the choice of $X_\zeta \subseteq   X$, and hence
$(*)_\zeta$ holds for some $\zeta<\lambda$.

\medskip

By induction on $n<\omega$ we now choose ordinals $\alpha_n, \beta_n, \gamma_n, \delta_n$ in $\delta$ and $ \sigma(n)<\lambda^+$ so that the following conditions are satisfied:
\begin{description}
\item{A:\ } $\langle  (\alpha_m, \beta_m, \gamma_m, \delta_m):m\leq n \rangle$ is an intial segment of a play in the game ${\cal   G}_\zeta$ in which $P_1$ uses the winning strategy $\tau_\zeta$.

\smallskip

\item{B:\ } $\alpha_0, \beta_0<\min(C_ \delta)$.

\smallskip

\item{C:\ } $\gamma_n=\min\{\gamma: \gamma>i_{\sigma(2n)} \wedge
\gamma\, e_{X_n}\,\delta \wedge c(\gamma,\delta)=\zeta\}$,
where 
$$X_n=  \bigcup\{\{\alpha_\ell, \beta_ \ell, \gamma_{ \ell}, \delta_{\ell}\}:\ell<n\}\cup\{\alpha_n, \beta_n\}
\cup
\bigcup\{\{i_{\sigma(\ell)},i_{\sigma(\ell)+1}\}:\ell<2n\}.$$

\item{D:\ } $\delta_n=i_{\sigma(2n+1)}$.

\smallskip

\item{E:\ } For $n>0$, $[\alpha_n, \beta_n] \subseteq  (\delta_{n-1}, i_{
\sigma(2n-1)+1})$.

\smallskip

\item{F:\ }  $i_{\sigma(n)}$ belongs to $A_0$ or $A_1$ according as $n$ is even or odd and $\sigma(n)+1<\sigma(n+1)$.
\end{description}

We have to prove that it is possible to choose the $\alpha_n$ etc., so that these conditions are satisfied.  Clearly (B) holds since, by (oo), the first moves by $P_1$ using the stategy $\tau_\zeta$ are $\alpha_0<\beta_0
<\rho_0$ and $\rho_0\leq \min(C_ \delta)\in K$. 
By  $(*)_\zeta$, there are $\sigma(0)<\lambda^+$ and $\gamma$ such that
$i_ {\sigma(0)}\in A_0$, $i_{\sigma(0)}<\gamma<i_{\sigma(0)+1}$, $\gamma \,e_ {X_0}\, \delta$, where $X_0=\{\alpha_0,\beta_0\}$ and $c(\gamma,\delta)=\zeta$; let $\gamma_0$ be the least such $\gamma$. Now let $\sigma(1)>\sigma(0)+1$ be minimal so that
$i_{\sigma(1)}\in A_1$, and put $\delta_0=i_{\sigma(1)}$.
Now suppose that $n>0$ and that the $\alpha_m, \beta_m, \gamma_m, \delta_m$, $\sigma(2m)$ and $\sigma(2m+1)$ have been suitably defined for all $m<n$.  Let $\rho\in K$ be minimal such that $\rho>\delta_{n-1}$. 
$P_1$ chooses $\alpha_n, \beta_n$ using the strategy $\tau_\zeta$ 
so that $\delta_{n-1}<\alpha_n<\beta_n<\rho$.
Since $\delta_{n-1}=i_{\sigma(2n-1)}\in A_1$,  it follows that $i_{\sigma(2n-1)+1}\in K$ and hence $\rho\leq i_{\sigma(2n-1)+1}$.  Now by $(*)_\zeta$, there are $\sigma(2n)$ and $\gamma$ so that
$i_{\sigma(2n)}\in A_0$, 

 $i_{\sigma(2n)}<\gamma<i_{\sigma(2n)+1}$, $\gamma\, e_ {X_n}\,
\delta$ (where $X_n$ is as described in (C)), and
$c(\gamma,\delta)=\zeta$; let $\gamma_{n}$ be the least such $\gamma$.
Note that, since $i_{\sigma(2n)}\in A_0$, $i_{\sigma(2n)+1}=\min(C_
\delta \,\backslash\, (i_{\sigma(2n)}+1))\in K$.  Finally, choose a
minimal ordinal $\sigma(2n+1)>\sigma(2n)+1$ so that
$\delta_n=i_{\sigma(2n+1)}\in A_1$.  This completes the definition of
the $\alpha_n$ etc., so that (A)-(F) hold.

\medskip
By (C) it follows that $c(\gamma_n,\delta)=\zeta$ for all $n<\omega$, and hence $c(\gamma_m,\gamma_n)=\zeta$ holds for all $m<n<\omega$ since $\gamma_m\in X_n$ and $\gamma_n\, e_{X_n}\, \delta$.
There is no edge of  $G_ \kappa$ from $\delta$ to $(\alpha_0,\beta_0)$ since $\beta_0<\min(C_ \delta)$.    Since 
 $\gamma_n\, e_{X_n}\, \delta$ and $\beta_0\in X_n$, it follows that
 $\beta_0<\min (C_{\gamma_n})$ also, and so there is no edge from $\gamma_n$ to $(\alpha_0,\beta_0)$ either.  
By the construction, for $0<m<\omega$, $i_{\sigma(2m-1)}<\alpha_m<\beta_m<i_{\sigma(2m-1)+1}$, and hence $C_ \delta\cap(\alpha_m,\beta_m)=\emptyset$.  Therefore,
for any $\xi\in (\alpha_m,\beta_m)$, $h^*_ \delta(\sup(\xi\cap C_ \delta))=h^*_ \delta(i_{\sigma(2m-1)})=1$ by (F), and so there is no edge of $G$ from $\delta$ to $(\alpha_m, \beta_m)$.  If $0<m<n<\omega$, then $\gamma_n\, e_{X_n}\, \delta$ and therefore
,
$$ \mathop{\rm  tp}(\alpha_m\cap C_{\gamma_n})=
\mathop{\rm  tp}(\alpha_m\cap C_{\delta})=\mathop{\rm  tp}(\beta_m\cap C_{\delta})=\mathop{\rm  tp}(\beta_m\cap C_{\gamma_n}).$$  Therefore, 
for any $\xi\in (\alpha_m,\beta_m)$, it follows that
$$h^*_{ \gamma_n}(\sup(\xi\cap C_{\gamma_n}))=h^*_{\gamma_n}(\sup(\alpha_m\cap C_{\gamma_n}))=h^*_ \delta(\sup(\alpha_m\cap C_{\delta}))=1$$
and so
there are no edges of $G$ from $\gamma_n$ to $(\alpha_m, \beta_m)$ either.

\medskip

Thus we have produced a play in the game ${\cal   G}_\zeta$ in which 
$P_1$ uses the strategy $\tau_\zeta$ but the second player $P_2$ wins!
This contradicts the assumption that $\sigma_\zeta$ is a winning strategy for the first player, and completes the proof.\ \  
\end{proof}


\begin{thebibliography}{99}

\bibitem[]{ER} P. Erd\H{o}s and R. Rado,
{A partition calculus in set theory}, {\it Bull. Amer. Math.
Soc.}\/ {\bf 62} (1956) 427-489. 


\bibitem[]{HK} A. Hajnal and P. Komjath, 
{Embedding graphs into colored graphs}, {\it Trans. Amer.
Math. Soc.}\/ {\bf 307} (1988),    
395--409; Corrigendum: {\bf 332} (1992), 475.

\bibitem[]{KM} P. Komjath and E.C. Milner, {On a conjecture of R\" odl and Voigt.\/}
J. Combin. Theory, Ser. B {\bf 61} (1994), 199-209.

\bibitem[]{RV} V. R\" odl and B. Voigt, {Monochromatic trees with respect to
edge partitions}, {\it J. Combin. Theory Ser. B}\/ {\bf 58} (1993), 291-298.

\bibitem[]{Sh:289} Saharon Shelah [Sh: 289],
\newblock {Consistency of positive partition theorems for graphs and models},
\newblock in: {\it {Set theory and its applications (Toronto, ON, 1987)}},
{\em {Lecture Notes in Mathematics}} {\bf 1401},
(J.~Steprans and S.~Watson, eds.), 
\newblock {Springer, Berlin-New York}, (1989) 167--193.

\bibitem[]{Sh:365}Saharon Shelah [Sh: 365],
 {There are Jonsson algebras in many inaccessible cardinals},
in: {{\it Cardinal Arithmetic}}, {\em {Oxford Logic
  Guides}} {\bf 29} chapter III, {Oxford University Press}, 1994.


\bibitem[]{Sh:413}Saharon Shelah [Sh: 413],
\newblock {More Jonsson Algebras and Colourings},
\newblock {\it {Archive for Mathematical Logic}}, to appear.

\bibitem[]{Sh:572}Saharon Shelah [Sh: 572],
\newblock {Colouring and $\aleph_2$-cc not productive},
\newblock {\it {Annals of Pure and Applied Logic}},
\newblock {{{\bf 84} (1997), 153-174.}}.
\end{thebibliography}
\end{document}